\documentclass[12pt,a4paper]{article}
\usepackage[T1]{fontenc}
\usepackage[utf8]{inputenc}
\usepackage{lmodern,amsmath,amssymb,amsthm,xcolor}
\usepackage[margin=25mm]{geometry}
\usepackage{microtype,enumitem}
\usepackage[hidelinks,hyperfootnotes=false]{hyperref}
\newtheorem{theorem}{Theorem}[section]
\newtheorem{conjecture}[theorem]{Conjecture}
\newtheorem{proposition}[theorem]{Proposition}
\newtheorem{lemma}[theorem]{Lemma}
\newtheorem{corollary}[theorem]{Corollary}
\theoremstyle{definition}
\newtheorem{example}[theorem]{Example}
\theoremstyle{plain}
\newtheorem{claim}{Claim}
\numberwithin{equation}{section}

\title{The size and structure of graphs with a unique $k$-factor\thanks{This work was supported by National Natural Science Foundation of China (Nos. 12401464 and 12471334).}}
\author{Zixuan Yang\\[3pt]
\small School of Mathematics, Northwest University, Xi'an, Shaanxi, P.R. China
}
\date{}

\begin{document}

\maketitle
{\renewcommand{\thefootnote}{}\footnotetext{E-mail addresses: \texttt{yangzixuan@nwpu.edu.cn} (Z. Yang)}}

\begin{abstract}
For integers $n>k\ge1$ with $nk$ even, let $m(n,k)$ be the maximum size of a graph of order $n$ with a unique $k$-factor. In 1964, Hetyei and in 1984, Hendry [\emph{J. Combin. Theory Ser. B} \textbf{37} (1984), 53--63] proved $m(n,1)=n^2/4$ for even $n$ and $m(n,2)=\lfloor n(n+1)/4\rfloor$, respectively. In 2000, Johann [\emph{J. Graph Theory} \textbf{35} (2000), 227--243] showed that $m(n,k)=nk/2+\binom{n-k}{2}$ for $k>n/2$ and $kn$ even and $m(n,k)=n^2/4+(k-1)n/4$ for $n=2kq$, where $q$ is a positive integer. Later, Volkmann [\emph{Combinatorica} \textbf{24} (2004), 531--540] proved $m(n,k)=k^2+\binom{n-k}{2}$ for $n/3\le k<n/2$ and $kn$ even, and he determined $m(n,3)$. In this paper, we determine $m(n,k)$ for all admissible pairs, proving the size conjecture posed by Volkmann. We also prove that every extremal graph with $2\le k\le n-2$ has a universal vertex whose neighbors in the unique factor are precisely its degree-$k$ vertices. In particular, this result confirms Volkmann's structural conjecture.

\end{abstract}

\section{Introduction}
Throughout this paper, we consider finite, undirected, simple graphs. The order and size of a graph $G$ are $|V(G)|$ and $|E(G)|$, respectively. For $x\in V(G)$, let $d_G(x)$ denote its degree. For $X\subseteq V(G)$, let $G[X]$ be the induced subgraph; for disjoint $X,Y\subseteq V(G)$, let $e_G(X,Y)$ count the edges between them. Write $e_G(X)=|E(G[X])|$. A vertex adjacent to every other vertex is universal.
We denote by $K_n$ the complete graph of order $n$, by $K_{p,q}$ the complete bipartite graph with partite sets of sizes $p$ and $q$, and by $\overline G$ the complement of $G$. For disjoint graphs $G$ and $H$, their join $G+H$ is obtained from their disjoint union by adding all edges between $V(G)$ and $V(H)$.

For a function $f:V(G)\to \mathbb Z_{\ge0}$, an $f$-factor is a spanning subgraph $F$ satisfying $d_F(x)=f(x)$ for every $x\in V(G)$. If $f(x)=k$ for every vertex $x\in V(G)$, then $F$ is called a $k$-factor. For integers $n>k\ge1$ with $nk$ even, define
\[
 m(n,k)=\max\{|E(G)|:|V(G)|=n\text{ and }G\text{ has a unique }k\text{-factor}\}.
\]
Whenever $F$ is the unique $k$-factor, we color its edges red and the remaining edges of $G$ blue, and write $E_r=E(F)$ and $E_b=E(G)\setminus E(F)$. We use $N_r(x)$ and $N_b(x)$ for the red and blue neighborhoods of the vertex $x\in V(G)$.

In 1964, Hetyei (see Corollary 1.6 in \cite{Lovasz1972}) determined $m(n,1)$, and in 1984, Hendry \cite{Hendry} determined $m(n,2)$. In 2000, Johann \cite{Johann} proved the general bound
\[
 m(n,k)\le\frac{n^2+(k-1)n}{4},
\]
with equality when $n$ is a multiple of $2k$, and determined the maximum graphs when $k>n/2$. Later, Volkmann \cite{Volkmann} determined $m(n,k)$ for $2k<n\le3k$ and for $k=3$, and supplied constructions for all suitable positive integers $n$ and $k$.

Let $n=2kq+r$, $q\ge0$, and $0\le r<2k$, and define
\[
 \vartheta_k(r):=\begin{cases}
 r(k+1-r),&0\le r\le k,\\
 (r-k-1)(2k-r),&k+1\le r<2k.
 \end{cases}
\]
Define
\[
 h(n,k):=\frac{n^2+(k-1)n-\vartheta_k(r)}4.
\]
Volkmann \cite{Volkmann} conjectured that $m(n,k)=h(n,k)$ for all suitable positive integers $n$ and $k$.
\begin{conjecture}[Volkmann, \cite{Volkmann}]\label{conj:size}
For positive integers $n>k\ge2$ and $nk$ even, $m(n,k)=h(n,k)$.
\end{conjecture}
Our first result proves Conjecture~\ref{conj:size}.
\begin{theorem}\label{thm:size}
Let $n>k\ge2$ be integers such that $nk$ is even.
If $n=2kq+r$ with $0\le r<2k$, then
\[
 m(n,k)=
 \begin{cases}
 \displaystyle
 \frac{n^2+(k-1)n-r(k+1-r)}4,
 &\text{for}~~0\le r\le k,\\[6pt]
 \displaystyle
 \frac{n^2+(k-1)n-(r-k-1)(2k-r)}4,
 &\text{for}~~k+1\le r<2k.
 \end{cases}
\]
\end{theorem}

\noindent\textit{Remark.} The previously determined cases include $k=2,3$, the range $n\le3k$, and the residue classes $r=0,k+1$. Outside these ranges, the unresolved cases of the conjecture are contained in $k\ge4$, $n>3k$, and $r\in\{1,\ldots,k\}\cup\{k+2,\ldots,2k-1\}$, subject to $nk$ being even.

Hoffmann and Volkmann \cite{HoffmannVolkmann2004,HoffmannVolkmann2005} studied unique regular factors and saturated graphs with prescribed degrees in bipartite graphs. 
Hoffmann, Sidorowicz, and Volkmann \cite{HoffmannSidorowiczVolkmann} solved the extremal size problem in bipartite graphs and determined structural properties of its extremal graphs. Volkmann~\cite[Corollary~2.3]{Volkmann} proved that a graph with a unique $k$-factor, $k\ge2$, has minimum degree at most $k+1$; see also Akbari, Ghodrati, and Hosseinzadeh~\cite[Corollary~1.5]{Akbari}.
For a graph with a unique $k$-factor, Volkmann \cite{Volkmann} gave the following conjecture.

\begin{conjecture}[Volkmann, \cite{Volkmann}]\label{conj:structure}
Let $G$ be a graph of order $n$ with a unique $k$-factor such that $2\le k\le n-2$. If $G$ is of maximum size, then $G$ contains exactly $k$ vertices of degree $k$.
\end{conjecture}
Our second result shows the following stronger statement.
\begin{theorem}\label{thm:structure}
Let $2\le k\le n-2$, and let $G$ be a graph of order $n$ and size $m(n,k)$ with unique $k$-factor $F$. Then $D=\{v\in V(G):d_G(v)=k\}$ has size $k$. Moreover, $G$ has a universal vertex $x$ such that $N_r(x)=D$.
\end{theorem}

\section{Some technical lemmas}
The following theorem will be used when every prescribed degree is positive.

\begin{theorem}[Jackson and Whitty, \cite{JacksonWhitty}]\label{thm:JW}
Let $G$ be a 2-edge-connected graph with a unique $f$-factor, where $f(v)>0$ for every $v\in V(G)$. Then some vertex $x$ satisfies $d_G(x)=f(x)$.
\end{theorem}

We recall the following necessary factor inequality, which is a corollary of Lov\'asz's $(g,f)$-factor theorem~\cite{Lovasz1970} with $g(v)=f(v)=b(v)$ for every vertex $v$.
 For a function $b$, write $b(X)=\sum_{x\in X}b(x)$.

\begin{proposition}\label{prop:factor}
Let $H$ be a graph with a $b$-factor, where $b:V(H)\to \mathbb Z_{\ge0}$. 
If $S,T$ partition $V(H)$, then
\[
 b(S)-b(T)+\sum_{v\in T}d_H(v)-e_H(S,T)\ge0.
\]
\end{proposition}

Fix a finite vertex set $U$ and $f:U\to \mathbb Z_{\ge0}$. We call a graph $R$ on $U$ \emph{$f$-extremal} if it has a unique $f$-factor and has maximum size among graphs on $U$ with a unique $f$-factor.  
For $X\subseteq U$, write $f(X)=\sum_{v\in X}f(v)$.
 If $L$ is its unique factor, a vertex is \emph{free} if $f(v)=0$, and \emph{tight} if $f(v)>0$ and $d_R(v)=f(v)$.

\begin{lemma}\label{lem:delete}
Let $R$ be $f$-extremal with unique $f$-factor $L$. Then the following
two statements hold:
\begin{enumerate}
\item[(i)] if $x$ is free, then $x$ is universal, and $R-x$ is extremal
for $f|_{U\setminus\{x\}}$;
\item[(ii)] if $x$ is tight, then every edge
incident to $x$ belongs to $L$, and $R-x$ is extremal
for the function
$f_x(y)=f(y)-\mathbf1_{xy\in E(L)}$
for every $y\in U\setminus\{x\}$.
\end{enumerate}
\end{lemma}

\begin{proof}
We prove the two statements separately.

\noindent (i) Let $x$ be a free vertex.
Since $f(x)=0$, every $f$-factor has degree $0$ at $x$.
It follows that no edge incident to $x$ can be used in any $f$-factor
of $R$.

We now show that $x$ is universal. Suppose, for contradiction, that
$x$ is not universal. Then there exists some vertex
$y\in U\setminus\{x\}$ with $xy\notin E(R)$.
Let $R'=R+xy$.
Since edges incident to $x$ cannot belong to any $f$-factor, the set
of $f$-factors of $R'$ is exactly the same as the set of $f$-factors
of $R$. In particular, $R'$ still has a unique $f$-factor $L$.
However,
$|E(R')|=|E(R)|+1>|E(R)|$,
which contradicts the assumption that $R$ is $f$-extremal.
Hence $x$ is universal.

We next prove that $R-x$ is extremal for $f|_{U\setminus\{x\}}$.
Suppose, for contradiction, that $R-x$ is not
extremal. Note that $R-x$ has a unique
$f|_{U\setminus\{x\}}$-factor, namely $L-x$,
since every such factor extends to an $f$-factor
of $R$ by adjoining $x$ as an isolated vertex.
Thus there exists a graph $R''$ on vertex set
$U\setminus\{x\}$ such that
$|E(R'')|>|E(R-x)|$
and $R''$ has a unique $f|_{U\setminus\{x\}}$-factor.
Construct a new graph $R'''$ on $U$ by adding the vertex $x$ to
$R''$ and joining $x$ to every vertex in $U\setminus\{x\}$.
Since $x$ remains free, the $f$-factors of $R'''$ are in bijection
with the $f|_{U\setminus\{x\}}$-factors of $R''$, so $R'''$
has a unique $f$-factor.
Since $x$ is universal in $R$, we have
\[
 |E(R''')|=|E(R'')|+(|U|-1)
 >|E(R-x)|+(|U|-1)=|E(R)|,
\]
contradicting the extremality of $R$.
Therefore $R-x$ must be $f|_{U\setminus\{x\}}$-extremal.

\noindent (ii) Let $x$ be a tight vertex. Then $f(x)>0$ and
$d_R(x)=f(x)$. Note that every $f$-factor of $R$ contains all
edges incident to $x$. Indeed, for any $f$-factor $F$, we have
$d_F(x)=f(x)=d_R(x)$.
As $F$ is a spanning subgraph of $R$, this equality
forces $F$ to contain all edges incident to $x$.

Also, $f_x(y)\ge0$ for every $y\ne x$, since
$xy\in E(L)$ implies $f(y)=d_L(y)\ge1$.
We now establish a bijection between the set of $f$-factors of $R$
and the set of $f_x$-factors of $R-x$. Let $F$ be any $f$-factor
of $R$. Since $F$ must contain all edges incident to $x$, the
subgraph $F-x$ satisfies
\[
 d_{F-x}(y)
 =d_F(y)-\mathbf1_{xy\in E(F)}
 =f(y)-\mathbf1_{xy\in E(L)}
 =f_x(y)
\]
for all $y\ne x$. Hence $F-x$ is an $f_x$-factor of $R-x$.
Conversely, let $F'$ be any $f_x$-factor of $R-x$.
Construct a spanning subgraph $F$ on $U$ by setting
$E(F)=E(F')\cup E_R(x)$, where $E_R(x)$ is the set of edges
incident to $x$ in $R$. Then $d_F(x)=d_R(x)=f(x)$ and
\[
 d_F(y)
 =d_{F'}(y)+\mathbf1_{xy\in E(L)}
 =f_x(y)+\mathbf1_{xy\in E(L)}
 =f(y)
\]
for all $y\ne x$, so $F$ is an $f$-factor of $R$.
These two constructions are inverse to each other.
Since $R$ has a unique $f$-factor, it follows immediately that
$R-x$ has a unique $f_x$-factor, namely $L-x$.

It remains to show that $R-x$ is $f_x$-extremal.
Suppose, for contradiction, that there exists a graph $R''$ on
$U\setminus\{x\}$ with
$|E(R'')|>|E(R-x)|$
and a unique $f_x$-factor.
Construct a new graph $R'''$ on $U$ by adding the
vertex $x$ to $R''$, and joining $x$ exactly to its neighbors in $L$.
Since $d_{R'''}(x)=f(x)$, every $f$-factor of $R'''$
contains all edges incident to $x$. Hence the same bijection
argument shows that
$R'''$ has a unique $f$-factor, and
\[
 |E(R''')|=|E(R'')|+f(x)
 >|E(R-x)|+f(x)=|E(R)|,
\]
which contradicts the extremality of $R$.
Hence $R-x$ is $f_x$-extremal.
\end{proof}

\begin{theorem}\label{thm:peel}
Let $R$ be an $f$-extremal graph of positive order with unique
$f$-factor $L$. If $f(v)>0$ for every $v\in V(R)$, then $R$ has
a tight vertex.
\end{theorem}

\begin{proof}
We use induction on $|V(R)|$. For $|V(R)|=1$, a graph of order one
cannot admit an $f$-factor with positive prescribed degree at its
only vertex, so the initial case is vacuous. Assume that the result
holds for all smaller orders.
In the following, we will complete the proof by contradiction. Suppose that $R$ has no tight vertex.

\begin{claim}\label{cl:connected}
$R$ is connected.
\end{claim}
\begin{proof}[Proof of Claim~\ref{cl:connected}]
Suppose to the contrary that
$R$ is disconnected, and let $A$ be the vertex set of one connected
component. Add an edge $e$ joining a vertex in $A$ to a vertex in
$V(R)\setminus A$, and denote the enlarged graph by $R'$.
Since $A$ is the vertex set of a component of $R$, the existing
factor $L$ gives $f(A)=2e_L(A)$, which is even.
For any $f$-factor $M$ of $R'$, the total degree sum over $A$
satisfies
\[
 f(A)=\sum_{v\in A}f(v)=2e_M(A)+\delta_M(A),
\]
where $\delta_M(A)$ is the number of edges in $M$ crossing the cut
between $A$ and $V(R)\setminus A$.
Since $2e_M(A)$ is even, $\delta_M(A)$ has the same parity as $f(A)$.
In $R'$, the only edge crossing this cut is the newly added edge $e$,
so $\delta_M(A)\in\{0,1\}$. Thus $e$ cannot be contained in any
$f$-factor of $R'$. It follows that the set of $f$-factors of $R'$
is identical to that of $R$, so $R'$ still has a unique $f$-factor.
However,
$|E(R')|=|E(R)|+1>|E(R)|$,
contradicting the assumption that $R$ is $f$-extremal.
Hence the claim holds.
\end{proof}

If $R$ has no bridge, then it is $2$-edge-connected, and by Theorem~\ref{thm:JW}, there exists a vertex $v$ with $d_R(v)=f(v)>0$, contradicting our assumption. So we assume that $R$ has a bridge.
Let $ac$ be a bridge of $R$, and let $A,C$ be the vertex sets of
the components of $R-\{ac\}$, with $a\in A$ and $c\in C$.
Let $\varepsilon=\mathbf1_{ac\in E(L)}$,
$f_A(v)=d_{L[A]}(v)$ for $v\in A$, and
$f_C(v)=d_{L[C]}(v)$ for $v\in C$.
Note that $f_A(a)=f(a)-\varepsilon$ and
$f_C(c)=f(c)-\varepsilon$, while $f_A(v)=f(v)$ for all
$v\in A\setminus\{a\}$ and $f_C(v)=f(v)$ for all
$v\in C\setminus\{c\}$.

For any $f$-factor $M$ of $R$, by the same parity argument as
before, we have
\begin{equation}\label{eq:bridge-parity}
 f(A)=2e_M(A)+\mathbf1_{ac\in E(M)}.
\end{equation}
Since the left-hand side is fixed,
$\mathbf1_{ac\in E(M)}$ is uniquely determined by parity.
It therefore equals $\varepsilon$, the value taken
by the unique factor $L$. In other words, every $f$-factor of $R$
contains the bridge $ac$ if and only if $L$ does, which establishes
a bijection between the $f$-factors of $R$ and pairs consisting of
an $f_A$-factor of $R[A]$ and an $f_C$-factor of $R[C]$, combined
with the fixed status of $ac$.
Since $R$ has a unique $f$-factor, both $R[A]$ and $R[C]$ have
unique factors for their respective degree functions%
, namely $L[A]$ and $L[C]$.
Moreover, both induced subgraphs are extremal.
Indeed, if either side could be replaced by a larger
graph with a unique factor for the same degree function, then $ac$
would still be the only edge between $A$ and $C$.
Its membership in an $f$-factor would still be fixed by the parity
identity~\eqref{eq:bridge-parity}. The same correspondence would
therefore yield a larger graph on $V(R)$ with a unique $f$-factor,
contradicting the extremality of $R$.

We now consider two cases based on whether the bridge belongs to
the unique factor.

\noindent\textbf{Case 1:} $ac\notin E(L)$.

Since $\varepsilon=0$, we have $f_A(a)=f(a)>0$, so $f_A(v)>0$
for all $v\in A$. By the induction hypothesis, the $f_A$-extremal
graph $R[A]$ contains at least one tight vertex.

By repeatedly applying Lemma~\ref{lem:delete}(ii), we may delete
tight vertices one by one from $R[A]$, while preserving the
extremal property for the residual degree function.
At each stage at which all residual prescribed degrees are
positive, the induction hypothesis supplies a tight vertex, since
the residual graph has smaller order than $R$.
We stop as soon as a free vertex appears.
Such a vertex must appear before all vertices are deleted:
every residual graph has its residual factor, and a factor on a
single vertex has degree zero.
Therefore, after deleting some sequence of tight vertices
$s_1,s_2,\dots,s_j$, we obtain a residual graph containing a free
vertex $x$.

We now construct a new graph $R^*$ from $R$ by removing the bridge
$ac$, and adding all edges between $x$ and vertices in $C$.
We claim that $R^*$ still has a unique $f$-factor.
Since $ac\notin E(L)$, the only removed edge lies outside $L$,
so $L$ remains an $f$-factor of $R^*$.
Since $x\notin\{s_1,\ldots,s_j\}$, no added edge is
incident to any $s_i$. Moreover, after removing $ac$ and the
preceding vertices $s_1,\ldots,s_{i-1}$, the edges incident to
$s_i$ are exactly those at the corresponding stage of the deletion
process in $R[A]$.
The vertex $s_1$ is therefore tight in $R^*$,
so all its incident edges are forced into any $f$-factor.
After deleting $s_1$, the same holds for $s_2$, and so on, up to
$s_j$. After removing all these tight vertices, the residual
prescribed degree of $x$ is $0$, so none of the new edges incident
to $x$ can be contained in any factor.
The residual factor on
$A\setminus\{s_1,\ldots,s_j\}$ is unique by the repeated
applications of Lemma~\ref{lem:delete}(ii), and the factor on $C$
is unique because $R[C]$ has a unique $f_C$-factor.
Hence the $f$-factor of $R^*$ is unique.

We now count edges. Note that we removed one edge (namely $ac$)
and added $|C|$ edges. So we have
\[
 |E(R^*)|=|E(R)|-1+|C|
 =|E(R)|+(|C|-1).
\]
Since $ac\notin E(L)$, the graph $L[C]$ is an
$f|_C$-factor of $R[C]$. As $f(v)>0$ for every $v\in C$,
it follows that $|C|\ge2$.
Thus $|E(R^*)|>|E(R)|$.
This contradicts the extremality of $R$.

\noindent\textbf{Case 2:} $ac\in E(L)$.

We first assert that $f(a)=1$ and $f(c)=1$.
Suppose for contradiction that $f(a)>1$.
Then $f_A(a)=f(a)-1\ge1$, so $f_A(v)>0$ for all $v\in A$.
By induction, $R[A]$ contains an $f_A$-tight vertex $v$.
If $v\ne a$, then
$d_R(v)=d_{R[A]}(v)=f_A(v)=f(v)$,
which implies that $v$ is tight in $R$, contradicting our
assumption. If $v=a$, then
$d_{R[A]}(a)=f_A(a)=f(a)-1$.
Consequently,
{
\[
 d_R(a)=d_{R[A]}(a)+1=f_A(a)+1=f(a),
\]
}
which implies that $a$ is tight in $R$, again a contradiction.
Hence $f(a)=1$, and by symmetry $f(c)=1$.

With $f(a)=f(c)=1$ and $\varepsilon=1$, we have $f_A(a)=0$ and
$f_C(c)=0$, which means that $a$ and $c$ are two free vertices
in $R[A]$ and $R[C]$, respectively.
By Lemma~\ref{lem:delete}(i), free vertices in extremal graphs are
universal. Therefore, $a$ and $c$ are adjacent to every vertex in
$A_0=A\setminus\{a\}$ and $C_0=C\setminus\{c\}$, respectively.

Notice that $A_0$ and $C_0$ are both nonempty.
Indeed, if $A_0=\emptyset$, then $a$ has only the
bridge $ac$ as its incident edge, which implies that
$d_R(a)=1=f(a)$.
Thus $a$ is a tight vertex, which contradicts our assumption.
The same holds for $C_0$.

We now construct a new graph $R^*$ as follows:
(i) retain all edges of $R[A_0]$ and $R[C_0]$,
and retain the bridge $ac$;
(ii) join $a$ to all vertices in $A_0\cup C_0$
(retaining the original edges from $\{a\}$ to $A_0$, and adding
edges from $a$ to every vertex in $C_0$);
(iii) remove all edges between $\{c\}$ and $C_0$, so that $c$ is
adjacent only to $a$; and
(iv) choose arbitrary vertices $x\in A_0$ and $y\in C_0$, and add
the edge $xy$.

We first verify the edge count.
In the original graph $R$, there are $|C_0|$ edges between
$\{c\}$ and $C_0$.
In $R^*$, we remove these $|C_0|$ edges, add $|C_0|$ edges from
$\{a\}$ to $C_0$, and add the edge $xy$.
All other edges are unchanged.
Since $ac$ was the only edge between $A$ and $C$, the added edges were not contained in $R$.
Thus, we have
\[
 |E(R^*)|=|E(R)|-|C_0|+|C_0|+1=|E(R)|+1,
\]
which means that $R^*$ has strictly more edges than $R$.

Since $ac\in E(L)$ and $f(c)=1$, none of the removed edges between
$\{c\}$ and $C_0$ belongs to $L$.
Thus $L$ remains an $f$-factor of $R^*$.
To contradict the extremality of $R$, it remains to show that
$R^*$ has a unique $f$-factor.
Since $d_{R^*}(c)=1$ and $f(c)=1$, the edge $ac$ must be included
in every $f$-factor.
Since $f(a)=1$, including $ac$ gives $a$ its prescribed degree.
Thus no other edge incident to $a$ can belong to any $f$-factor.
After excluding all other edges incident to $a$, the only edge
still available to a factor across the cut between $A_0$ and
$V(R^*)\setminus A_0$ is $xy$.
The sum of the prescribed degrees on $A_0$ is
$f(A_0)=2e_L(A_0)$, which is even.
Any factor must use an even number of edges across this cut,
so the single edge $xy$ cannot be contained in the factor.
Consequently, apart from $ac$, every edge of an $f$-factor belongs to $R[A_0]$ or $R[C_0]$. By Lemma~\ref{lem:delete}(i) applied to the free vertices $a$ and $c$ of $R[A]$ and $R[C]$, respectively, each of these graphs has a unique factor for the restriction of $f$.
Together with the forced edge $ac$, these two factors determine
the unique $f$-factor of $R^*$, contradicting the extremality of $R$.

The proof of Theorem~\ref{thm:peel} is complete.
\end{proof}

If a vertex of prescribed degree zero is present, then use Lemma~\ref{lem:delete}(i); otherwise use Theorem~\ref{thm:peel} and Lemma~\ref{lem:delete}(ii). Repeating this process until the vertex set is exhausted, we obtain the following corollary.

\begin{corollary}\label{cor:order}
Let $R$ be $f$-extremal with unique factor $L$. Its vertices admit an ordering $v_1,\dots,v_{|V(R)|}$ such that, with $R_i=R[\{v_i,\dots,v_{|V(R)|}\}]$, $L_i=L[\{v_i,\dots,v_{|V(R)|}\}]$, and $f_i(v)=d_{L_i}(v)$, each $R_i$ is $f_i$-extremal, and $v_i$ is either free and universal in $R_i$, or tight in $R_i$.
\end{corollary}

\section{Proof of Theorem~\ref{thm:size}}
We first prove an inequality for a bipartite graph whose vertices in one partite set are ordered.

\begin{lemma}\label{lem:incidence}
Let $\Gamma$ be a bipartite graph with partite sets $U$ and
$Q=\{1,\ldots,z\}$ such that $d_\Gamma(v)=k\ge1$ for each $v\in Q$.
If $|U|=t$ and $d_1\le\cdots\le d_t$ is the nondecreasing degree
sequence on $U$, then
\[
 \sum_{u\in U}c(u)\ge\sum_{i=1}^t d_i+
 \sum_{i=1}^t\left\lfloor\frac{t-i}{k}\right\rfloor d_i,
\]
where $c(u)=\max\bigl(N_\Gamma(u)\cup\{0\}\bigr)$ for $u\in U$.
\end{lemma}

\begin{proof}
If $z=0$, then both sides are zero.
So we may assume that $z\ge1$.
We form the bipartite graph $H$ with two partite sets $U$ and $Q$,
and with $uv\in E(H)$ if and only if $1\le v\le c(u)$.
Note that $\Gamma\subseteq H$.

\begin{claim}\label{cl:inc}
For every $X\subseteq U$, $\sum_{u\in X}d_\Gamma(u)\le\sum_{v=1}^z\min\{k,|N_H(v)\cap X|\}.$
\end{claim}

\begin{proof}[Proof of Claim~\ref{cl:inc}]
{
For every $v\in Q$, the inclusion $\Gamma\subseteq H$ and the
equality $d_\Gamma(v)=k$ show that
\[
 |N_\Gamma(v)\cap X|
 \le\min\{k,|N_H(v)\cap X|\}.
\]
Counting the edges of $\Gamma$ between $X$ and $Q$ from each
partite set, we obtain that
\[
 \sum_{u\in X}d_\Gamma(u)
 =\sum_{v=1}^z|N_\Gamma(v)\cap X|
 \le\sum_{v=1}^z\min\{k,|N_H(v)\cap X|\}.
\]
This proves the claim.
}
\end{proof}

Let $u_1,\ldots,u_t$ be an ordering of $U$ so that
$c(u_1)\ge\cdots\ge c(u_t)$.
For $v\in Q$, let
$N_v=|N_H(v)|=|\{u\in U:c(u)\ge v\}|$.
Then $k\le N_v\le t$.
Indeed, since $d_\Gamma(v)=k$, the vertex $v$ has $k$ neighbors
$u$ in $\Gamma$, and each such neighbor satisfies $c(u)\ge v$.
Moreover, we have $N_H(v)=\{u_1,\ldots,u_{N_v}\}$, and
\begin{equation}\label{eq:3.6}
 \sum_{u\in U}c(u)
 =\sum_{u\in U}\sum_{v=1}^{c(u)}1
 =\sum_{v=1}^z\sum_{u:c(u)\ge v}1
 =\sum_{v=1}^zN_v.
\end{equation}

Write $r_+=\max\{r,0\}$.
For every integer $h\ge1$, define
$
 X_h=\{u_i:hk<i\le t\}.
$
Thus $|X_h|=(t-hk)_+$.
When $hk<t$, this is a set of $t-hk$ vertices with the smallest
values of $c(u)$; when $hk\ge t$, it is empty.

We now apply the inequality in Claim~\ref{cl:inc} to a family of
subsets $X_h$ of $U$.
Since $N_H(v)=\{u_1,\ldots,u_{N_v}\}$, its
intersection with $X_h$ consists of those $u_i$ with
$hk<i\le N_v$.
Therefore
$
 |N_H(v)\cap X_h|=(N_v-hk)_+.
$
Note that every subset of $U$ of size $(t-hk)_+$ has degree
sum at least the sum of the $(t-hk)_+$ smallest degrees.
Hence
\[
 \sum_{u\in X_h}d_\Gamma(u)
 \ge\sum_{i=1}^{(t-hk)_+}d_i.
\]
Substituting $X=X_h$ into the inequality in Claim~\ref{cl:inc}
and combining with the above identities, we obtain that
\begin{equation}\label{eq:3.7}
 \sum_{i=1}^{(t-hk)_+}d_i
 \le\sum_{v=1}^z\min\{k,(N_v-hk)_+\}.
\end{equation}

For every integer $N\ge k$, we have
\[
 \sum_{h\ge1}\min\{k,(N-hk)_+\}=N-k.
\]
Indeed, write $N=ak+b$ with integers $a\ge1$ and $0\le b<k$.
If $1\le h\le a-1$, then $N-hk\ge k$, so
$\min\{k,N-hk\}=k$; this contributes $a-1$ terms equal to $k$.
If $h=a$, then $N-ak=b$, which implies that $\min\{k,b\}=b$.
If $h>a$, then $N-hk<0$, so $(N-hk)_+=0$ and the summand vanishes.
Summing these contributions, we obtain
\[
 (a-1)k+b=ak+b-k=N-k.
\]

For $hk\ge t$, both sides of~\eqref{eq:3.7} are zero.
Thus we may sum~\eqref{eq:3.7} over all $h\ge1$;
all sums involved have only finitely many nonzero terms.
On the right-hand side we interchange the order of summation
and apply the identity for each $v$,
\[
 \sum_{h\ge1}\sum_{v=1}^z\min\{k,(N_v-hk)_+\}
 =\sum_{v=1}^z\sum_{h\ge1}\min\{k,(N_v-hk)_+\}
 =\sum_{v=1}^z(N_v-k).
\]
Using \eqref{eq:3.6}, this simplifies to
\begin{equation}\label{eq:3.8}
 \sum_{v=1}^zN_v-kz=\sum_{u\in U}c(u)-kz.
\end{equation}

The sum of the left-hand side of \eqref{eq:3.7} over $h\ge1$ is
$\sum_{h\ge1}\sum_{i=1}^{(t-hk)_+}d_i$.
For a fixed index $i$, the term $d_i$ appears in the inner sum
if and only if $i\le t-hk$, or equivalently $h\le(t-i)/k$.
Since $h$ runs over positive integers, the number of such indices
$h$ is exactly $\lfloor(t-i)/k\rfloor$.
Interchanging the order of summation, we have
\begin{equation}\label{eq:3.9}
 \sum_{h\ge1}\sum_{i=1}^{(t-hk)_+}d_i
 =\sum_{i=1}^t\left\lfloor\frac{t-i}{k}\right\rfloor d_i.
\end{equation}
Combining \eqref{eq:3.7}, \eqref{eq:3.8}, and \eqref{eq:3.9},
we have
\begin{equation}\label{eq:3.10}
 \begin{aligned}
 \sum_{i=1}^t\left\lfloor\frac{t-i}{k}\right\rfloor d_i
 &=\sum_{h\ge1}\sum_{i=1}^{(t-hk)_+}d_i\\
 &\le\sum_{h\ge1}\sum_{v=1}^z\min\{k,(N_v-hk)_+\}\\
 &=\sum_{v=1}^zN_v-kz
 =\sum_{u\in U}c(u)-kz.
 \end{aligned}
\end{equation}
By the handshaking lemma for bipartite graphs, the total degree
sum over $Q$ equals the total degree sum over $U$, i.e.,
\[
 kz=\sum_{v\in Q}d_\Gamma(v)
 =\sum_{u\in U}d_\Gamma(u)
 =\sum_{i=1}^t d_i.
\]
Substituting the above equality into \eqref{eq:3.10}, we have
\[
 \sum_{u\in U}c(u)\ge\sum_{i=1}^t d_i+
 \sum_{i=1}^t\left\lfloor\frac{t-i}{k}\right\rfloor d_i.
\]
This completes the proof of Lemma~\ref{lem:incidence}.
\end{proof}

Let $k$ and $n$ be two positive integers such that $kn$ is even and $2k+1\le n\le3k$. To establish a lower bound for $m(n,k)$, we introduce some families of graphs of order $n$ with a unique $k$-factor.
\begin{example}[see Example 3.1 in \cite{Volkmann}]\label{ex:Q}
Let $k$ and $r$ be positive integers. We construct the graph $Q(n,k)$ that is of size $k^2+\binom{n-k}{2}$ and has a unique $k$-factor by considering two cases.
In both cases, the prescribed degrees on $B$ lie between $0$ and $k$ and sum to $k^2$. To realize them, list the $k$ vertices of $A$ cyclically $k$ times, partition the list into consecutive segments of the prescribed lengths, and join each vertex of $B$ to its segment. Each segment has distinct entries, and every vertex of $A$ has degree $k$.

\emph{Case 1.} $n=2k+2r$ and $2r\le k$.
Let $H$ be the bipartite graph with the partite sets $A$ and $B$ such that $|A|=k$ and $|B|=k+2r$ with
\begin{align*}
 B={}&\{u_1,u_2,\ldots,u_{2r-1}\}\cup\{v_1,v_2,\ldots,v_{k-2r+1}\}\\
 &\cup\{w_1,w_2,\ldots,w_r\}\cup\{y_1,y_2,\ldots,y_r\}.
\end{align*}
Next we connect the vertices of $A$ and $B$ by edges such that $d_H(x)=k$ for all $x\in A$, $d_H(u_i)=k$ for $1\le i\le2r-1$, $d_H(v_i)=k-r$ for $1\le i\le k-2r+1$, $d_H(w_i)=k-2r+1$ for $1\le i\le r$, and $d_H(y_i)=0$ for $1\le i\le r$. It is a simple matter to see that one can realize such a bipartite graph and
\[
 \sum_{x\in A}d_H(x)=\sum_{x\in B}d_H(x)=k^2.
\]
Now the graph $Q(n,k)$ of order $n=2k+2r$ originates from $H$ by connecting all pairs of vertices in $B$ by an edge.
Every $k$-factor contains all edges of $H$, since every vertex of $A$ has degree $k$. The vertices $u_i$ therefore have degree zero in the remaining factor on $B$. Each $y_i$ must be adjacent in this factor to all other vertices of $B\setminus\{u_1,\ldots,u_{2r-1}\}$, which has order $k+1$. This gives every $v_i$ its remaining degree $r$, and forces the vertices $w_i$ to induce $K_r$. Thus the remaining factor, and hence the $k$-factor of $Q(n,k)$, exists and is unique.

\emph{Case 2.} $n=2k+2r-1$ with $2r\le k+1$ and $k=2q$ even.
Let $H$ be the bipartite graph with the partite sets $A$ and $B$ such that $|A|=k$ and $|B|=k+2r-1$ with
\[
 B=\{u_1,u_2,\ldots,u_{q+r-1}\}\cup\{v_1,v_2,\ldots,v_{q-r+1}\}\cup\{w_1,w_2,\ldots,w_{2r-1}\}.
\]
Next we connect the vertices of $A$ and $B$ by edges such that $d_H(x)=k$ for all $x\in A$, $d_H(u_i)=k$ for $1\le i\le q+r-1$, $d_H(v_i)=k-2r+1$ for $1\le i\le q-r+1$, and $d_H(w_i)=q-r+1$ for $1\le i\le2r-1$. Clearly, one can realize such a bipartite graph and
\[
 \sum_{x\in A}d_H(x)=\sum_{x\in B}d_H(x)=k^2.
\]
Now the graph $Q(n,k)$ of order $n=2k+2r-1$ originates from $H$ by connecting all pairs of vertices in $B$ by an edge.
Again, every $k$-factor contains all edges of $H$, and the vertices $u_i$ have degree zero in the remaining factor on $B$. Each $w_i$ must be adjacent in this factor to all other vertices of $B\setminus\{u_1,\ldots,u_{q+r-1}\}$, which has order $q+r$. These edges give every $v_i$ its remaining degree $2r-1$. Thus the remaining factor exists and is unique in this case as well.
\end{example}
We next introduce another family of graphs $H(n,k)$.
\begin{example}[see Example 5.1 in \cite{Volkmann}]\label{ex:H}
For $k<n\le2k$ with $nk$ even, if $P$ is a $(2k-n)$-regular graph of order $k$, then the graph $R(n,k)=P+K_{n-k}$ is of size $nk/2+\binom{n-k}{2}$ and has a unique $k$-factor. 
Such a graph $P$ exists since $0\le2k-n\le k-1$ and $k(2k-n)=2k^2-nk$ is even.
 Now let $n,k$, and $q$ be positive integers such that $n>2k$ and $n=2kq+r$ with $0\le r<2k$ and $nk$ even. Construct different families of graphs $H(n,k)$ of order $n$ as follows: take $q-1$ disjoint copies of $K_k+\overline{K_k}$ labelled $H_1,H_2,\ldots,H_{q-1}$.

\emph{Case 1.} For $r=0$ define $H_q=K_k+\overline{K_k}$.

\emph{Case 2.} For $1\le r\le k$ define $H_q=Q(2k+r,k)$ (cf. Example~\ref{ex:Q}).

\emph{Case 3.} For $k+1\le r<2k$ define $H_q=K_k+\overline{K_k}$ and $H_{q+1}=R(r,k)$.

For each $i$ and $j$ with $1\le i<j\le q$ in the first two cases and $1\le i<j\le q+1$ in the third case, add the edges connecting all vertices of the $K_k$ in $H_i$ to all the vertices of $H_j$.

It can easily be verified that $H(n,k)$ has a unique $k$-factor, and its size is $h(n,k)$. For $n>2k$, the construction uses ordered blocks. A standard block has two sets $A_i,B_i$ of size $k$, with $A_i$ independent, $B_i$ complete, and all edges between them. Its unique $k$-factor is $K_{k,k}$. Add all edges from $B_i$ to every later block.

When $r=0$, use $q$ standard blocks. When $1\le r\le k$, use $q-1$ standard blocks followed by the graph $Q(2k+r,k)$. When $k+1\le r<2k$, use $q$ standard blocks followed by Hendry's graph $R(r,k)$. In the range $k<n\le2k$, use $R(n,k)$ itself, with $R(2k,k)$ understood as a standard block.

The degree-$k$ vertices in the first standard block force its $K_{k,k}$ factor, giving every vertex of that block degree $k$ in the factor. No edge from that block to a later block can be used. Repetition reduces uniqueness to the last block. Counting the edges gives $h(n,k)$, as in Example~\ref{ex:H}.
\begin{samepage}
Indeed, prepending a standard block to a graph of order $N-2k$ adds
\[
 k^2+\binom{k}{2}+k(N-2k)
 =kN-\frac{k(k+1)}2
 =h(N,k)-h(N-2k,k)
\]
edges; the final blocks have the stated sizes.\end{samepage}
In particular,
\begin{equation}\label{eq:lower}
 m(n,k)\ge h(n,k).
\end{equation}
\end{example}

\begin{proof}[Proof of Theorem~\ref{thm:size}]
Let $G$ be a graph of order $n$ and size $m(n,k)$ with unique
$k$-factor $F$.
Color the edges of $F$ red and all other edges of
$G$ blue.
Write $n=2kq+r$, where $q\ge0$ and $0\le r<2k$.
Use the ordering $v_1,\ldots,v_n$ from Corollary~\ref{cor:order}.
Let $a_i=d_{F[\{v_i,\ldots,v_n\}]}(v_i)$.
Call the deletion of $v_i$ free when $a_i=0$ and tight otherwise.
At each step, a free vertex is universal in the
remaining graph, while a tight vertex has remaining degree $a_i$.
Let $s_1,\ldots,s_z$ be the free vertices in deletion order, with
positions $p_1<\cdots<p_z$, and let $T$ be the set of the other
$t=n-z$ vertices.

If $i<j$ and $v_i$ is free, then $v_iv_j$ is blue; and if $v_i$
is tight, every remaining edge incident to $v_i$
is red. Consequently, with $B=|E_b|$, we have
\begin{equation}\label{eq:blue-count}
 B=\sum_{i:a_i=0}(n-i)=\sum_{j=1}^z(n-p_j).
\end{equation}
Note that the free vertices are independent in $F$.
Indeed, the earlier endpoint of a red edge between
two free vertices would have positive residual factor degree when
deleted, contradicting the definition of a free deletion.

Let $X_j=N_r(\{s_1,\ldots,s_j\})\subseteq T$, $x_j=|X_j|$, and
$C=\sum_{j=1}^z x_j$.
Note that $p_j\ge j+x_j$,
since every red neighbor of $s_h$ precedes $s_h$,
and hence every vertex of $X_j$ precedes $s_j$.
So by \eqref{eq:blue-count}, we have
\begin{equation}\label{eq:blue-bound}
 B=\sum_{j=1}^z(n-p_j)
 \le\sum_{j=1}^z(n-j-x_j)
 =zn-\binom{z+1}{2}-C.
\end{equation}

We now construct an auxiliary bipartite graph to apply
Lemma~\ref{lem:incidence}.
For $u\in T$, let $d(u)$ be the number of red neighbors among the
free vertices. Since each free vertex has $k$ red neighbors, we
have $\sum_{u\in T}d(u)=kz$.
Define a bipartite graph $\Gamma$ with partite sets $T$ and
$\{1,\ldots,z\}$ by joining $u\in T$ to $v\in\{1,\ldots,z\}$
when $u$ is a red neighbor of $s_{z-v+1}$.
Every vertex in the second partite set has degree $k$.
As in Lemma~\ref{lem:incidence}, set
$c(u)=\max(N_\Gamma(u)\cup\{0\})$.
Double-counting the pairs $(u,j)$ with $u\in X_j$ shows
$
 C=\sum_{u\in T}c(u),
$
since, if $s_j$ is the first free red neighbor of $u$, then $u$
belongs to exactly $X_j,\ldots,X_z$ and $c(u)=z-j+1$;
if $u$ has no free red neighbor, both contributions are zero.

Let $d_1\le\cdots\le d_t$ be
the degree sequence of $\Gamma$ on $T$, and let
$\rho_i=k-d_i$, $w_i=\lfloor(t-i)/k\rfloor$,
$\ell=t-z=n-2z$, and $\Phi=\sum_{i=1}^t w_i\rho_i$.
Then $\rho_1\ge\cdots\ge\rho_t$ is the non-increasing degree
sequence of $F[T]$,
\[
 \sum_{i=1}^t\rho_i=kt-kz=k(n-2z)=k\ell,
\]
and $0\le\ell\le t$
(since combining $d(u)\le k$ and
$\sum_{u\in T}d(u)=kz$ gives $z\le t$).
By Lemma~\ref{lem:incidence}, we have
\begin{equation}\label{eq:C-rho}
 C\ge kz+k\sum_{i=1}^t w_i-\Phi.
\end{equation}
For $1\le j\le t$, write $D_j=\sum_{i=1}^j\rho_i$.
Then by the same counting argument
as in Lemma~\ref{lem:incidence}, we have
\begin{equation}\label{eq:Phi-partials}
 \Phi=\sum_{h=1}^{\lfloor(t-1)/k\rfloor}D_{t-hk}.
\end{equation}
Let $Y$ be a set of $j$ vertices of largest degree
in $F[T]$.
Write $D_j=\sum_{u\in Y}d_{F[T]}(u)$.
We decompose the degree sum over $Y$ as
\[
 D_j=2e_{F[T]}(Y)+e_{F[T]}(Y,T\setminus Y).
\]
The degree sum over the complement is
\[
 k\ell-D_j
 =2e_{F[T]}(T\setminus Y)+e_{F[T]}(Y,T\setminus Y),
\]
which implies that
$k\ell-D_j\ge e_{F[T]}(Y,T\setminus Y)$.
Since $F[T]$ is simple, $2e_{F[T]}(Y)\le j(j-1)$.
Therefore
\[
 D_j\le j(j-1)+e_{F[T]}(Y,T\setminus Y)
 \le j(j-1)+(k\ell-D_j),
\]
i.e., $2D_j\le k\ell+j(j-1)$.
Moreover, each vertex in $F[T]$ has degree at most $k$, so
$D_j\le kj$, and the total degree sum of $F[T]$ is $k\ell$,
so $D_j\le k\ell$. Thus
\begin{equation}\label{eq:partial-bound}
 D_j\le
 \min\left\{kj,\ k\ell,\ \frac{k\ell+j(j-1)}2\right\}.
\end{equation}

To bound $\Phi$, we first use only the numerical
constraints $0\le\rho_i\le k$ and $\sum_i\rho_i=k\ell$,
together with the fact that the weights $w_i$ are non-increasing.
Under these numerical constraints, if $i<j$, $\rho_i<k$, and
$\rho_j>0$, replace $\rho_i$ and $\rho_j$ by $\rho_i+\delta$
and $\rho_j-\delta$, respectively, where
$\delta=\min\{k-\rho_i,\rho_j\}$.
This preserves the constraints and does not decrease the weighted
sum. Repeating these exchanges, we obtain $k$ in each of the first
$\ell$ positions and $0$ in all remaining positions.
Hence $\Phi\le k\sum_{i=1}^{\ell}w_i$.
By \eqref{eq:C-rho}, we have
\begin{equation}\label{eq:C0}
 C\ge kz+k\sum_{i=\ell+1}^t w_i
 =k\sum_{j=1}^z\left\lceil\frac jk\right\rceil
 =:C_0(z,k).
\end{equation}
Here the equality uses $t-\ell=z$ and
$1+\lfloor(j-1)/k\rfloor=\lceil j/k\rceil$.
Combining \eqref{eq:blue-bound} and \eqref{eq:C0}, we have
\begin{equation}\label{eq:U}
 B\le U_{n,k}(z)
 :=zn-\binom{z+1}{2}-C_0(z,k).
\end{equation}
Then the first difference is
\begin{equation}\label{eq:U-difference}
 U_{n,k}(z+1)-U_{n,k}(z)
 =n-z-1-k\left\lceil\frac{z+1}{k}\right\rceil,
\end{equation}
which is strictly decreasing in $z$.

Let $J=n(n-k-1)/4$, and let
$\beta=h(n,k)-nk/2=J-\vartheta_k(r)/4$.
If $z=ak+b$, with $0\le b<k$, then
\[
 C_0(z,k)=\frac{k^2a(a+1)}2+kb(a+1).
\]
Using $n=2z+\ell$ and the expansion of \eqref{eq:U},
we obtain that
\begin{equation}\label{eq:J-U}
 J-U_{n,k}(z)
 =\frac{\ell(\ell-k-1)+2b(k-b)}4.
\end{equation}

Suppose first that $k<n\le2k$, and let $h_0=n-k$.
If $z=0$, then $B=0$; and if $z\ge1$, then the $k$ red neighbors
of each free vertex belong to $T$, giving
$z\le n-k=h_0\le k$.
In this range $\lceil j/k\rceil=1$
$(1\le j\le z)$, and hence
$C_0(z,k)=kz$.
So
\[
 B\le zn-\binom{z+1}{2}-kz
 =zh_0-\binom{z+1}{2}
 =\binom{h_0}{2}-\binom{h_0-z}{2}
 \le\binom{n-k}{2}.
\]
It follows that
$|E(G)|\le nk/2+\binom{n-k}{2}=h(n,k)$,
including $n=2k$.

Now assume that $n>2k$.
If $r=0$, then we let $z_0=kq$.
By \eqref{eq:U-difference},
\[
 U_{n,k}(z_0)-U_{n,k}(z_0-1)
 =n-z_0-k\left\lceil\frac{z_0}{k}\right\rceil
 =2kq-kq-kq=0,
\]
which, since the first differences are strictly decreasing, shows
that $z_0-1$ and $z_0$ maximize $U_{n,k}$.
At $z_0$, we have $\ell=b=0$.
So by \eqref{eq:J-U}, we obtain that
$B\le U_{n,k}(z_0)=J=\beta$.

If $k<r<2k$, then we let $z_0=kq+r-k$.
Since $\lceil z_0/k\rceil=q+1$,
by \eqref{eq:U-difference}, we have
{
\[
 U_{n,k}(z_0)-U_{n,k}(z_0-1)
 =n-z_0-k(q+1)=0.
\]
}
Thus $z_0-1$ and $z_0$ again maximize $U_{n,k}$.
At $z_0$, we have $\ell=2k-r$ and $b=r-k$,
whence, by \eqref{eq:J-U},
\[
 J-U_{n,k}(z_0)
 =\frac{(r-k-1)(2k-r)}4.
\]
Therefore $B\le\beta$ also holds.

It remains to consider $1\le r\le k$.
We discuss two cases according to $\ell$.

\medskip\noindent\textbf{Case 1: $\ell\ge k+1$.}

Since $k\lceil(z+1)/k\rceil\le z+k$,
by \eqref{eq:U-difference}, we have
$U_{n,k}(z+1)-U_{n,k}(z)\ge\ell-k-1$.
Thus $U_{n,k}(z)$ is nondecreasing on the integer range
$n-2z\ge k+1$.
Let $\ell_0$ be the least integer at least $k+1$ of the same parity
as $n$.
Let $\ell_0=k+\varepsilon$, where
\[
 \varepsilon=
 \begin{cases}
 1,&\text{if }r\not\equiv k\pmod2,\\
 2,&\text{if }r\equiv k\pmod2.
 \end{cases}
\]
The largest $z$ in this range is
$z_0=(n-\ell_0)/2=k(q-1)+b$, and
$b=(r+k-\varepsilon)/2$, where $0\le b<k$.
Subtracting $\vartheta_k(r)/4=r(k+1-r)/4$
from \eqref{eq:J-U}, we have
\[
 \beta-U_{n,k}(z_0)=
 \begin{cases}
 \dfrac{(k-r)^2-1}{8},
 &\text{if }\varepsilon=1,\\[6pt]
 \dfrac{(k-r)^2+2(k+r)}8,
 &\text{if }\varepsilon=2.
 \end{cases}
\]
When $\varepsilon=1$, $k-r$ is a positive odd integer, so the first
expression is nonnegative.
When $\varepsilon=2$, the second expression is positive.
Hence
$B\le U_{n,k}(z)\le U_{n,k}(z_0)\le\beta$.

\medskip\noindent\textbf{Case 2: $0\le\ell\le k$.}

Write $t=ak+s$, where $a\ge1$ and $1\le s\le k$;
if $k$ divides $t$, take $s=k$.
Here $a\ge1$ since $t=(n+\ell)/2>k$.
In particular, $a=\lfloor(t-1)/k\rfloor$, and \eqref{eq:Phi-partials} can be rewritten as  
$\Phi=\sum_{j=0}^{a-1}D_{jk+s}$.
In \eqref{eq:Phi-partials}, all the partial sums except the last
have indices at least $k+1$, and the last is $D_s$.
Apply $D_j\le k\ell$ to the first $a-1$ partial sums,
and apply \eqref{eq:partial-bound} to $D_s$,
to obtain that
\begin{equation}\label{eq:Phi-case2}
 \Phi\le(a-1)k\ell+\frac{k\ell+s(s-1)}2.
\end{equation}
Also,
\begin{equation}\label{eq:weights-sum}
 \sum_{i=1}^t w_i
 =\sum_{m=0}^{ak+s-1}\left\lfloor\frac mk\right\rfloor
 =k(0+1+\cdots+(a-1))+as
 =as+\frac{ka(a-1)}2.
\end{equation}
Substituting \eqref{eq:Phi-case2} and \eqref{eq:weights-sum}
into \eqref{eq:C-rho} and then \eqref{eq:blue-bound},
and using $z=ak+s-\ell$, we obtain that
\begin{align*}
 4(J-B)
 &\ge
 \ell(\ell-k-1)-2z(z+k)+4C\\
 &\ge\ell(\ell-k-1)-2z(z+k)+4kz
 +4k\left(as+\frac{ka(a-1)}2\right)\\
 &\qquad-4(a-1)k\ell-2k\ell-2s(s-1)\\
 &=(2s-\ell)(k+1-2s+\ell).
\end{align*}
Since $n=2z+\ell=2ak+2s-\ell$ and $n=2kq+r$,
we have $2s-\ell\equiv r\pmod{2k}$.
On the other hand, since $1\le s\le k$ and $0\le\ell\le k$,
$2-k\le2s-\ell\le2k$.
For $1\le r\le k$, we have
$r-2k\le-k<2-k$ and $r+2k>2k$.
Thus the only integer in this interval congruent to $r$ modulo
$2k$ is $r$ itself, so $2s-\ell=r$.
Therefore
$
 J-B\ge (r(k+1-r))/4,
$
which is exactly $B\le\beta$.

Since the unique factor has $nk/2$ edges,
\[
 m(n,k)=|E(G)|=\frac{nk}2+B\le h(n,k).
\]
The examples in this section give $m(n,k)\ge h(n,k)$.
Hence $m(n,k)=h(n,k)$. The proof of Theorem~\ref{thm:size} is complete.
\end{proof}

\section{Proof of Theorem~\ref{thm:structure}}
The proof of Theorem~\ref{thm:structure} relies on the deletion ordering from Corollary~\ref{cor:order} and the extremal deletion properties from Lemma~\ref{lem:delete}.

\begin{proof}[Proof of Theorem~\ref{thm:structure}]
Starting from $G$ with the unique $k$-factor $F$, we repeatedly
delete tight vertices following the deletion ordering from
Corollary~\ref{cor:order},
stopping at the first free deletion.
Such a step must occur, since the residual factor on the last
remaining vertex has degree zero.
Let $t_1,\ldots,t_h$ be the vertices deleted before this step,
listed in deletion order, and let $T_0=\{t_1,\ldots,t_h\}$.
Let $x$ be the vertex deleted at this first free step.

By Lemma~\ref{lem:delete}(ii), every edge incident to a tight
vertex in the corresponding residual graph belongs to the
residual factor (i.e., is a red edge).
For every edge of $G$ incident to a vertex of
$T_0$, its earlier endpoint in the deletion ordering also belongs
to $T_0$. The edge is still present when that endpoint is
deleted tightly, so it is red.
As a result, every vertex in $T_0$ has degree exactly $k$ in
the original graph $G$.

\setcounter{claim}{2}
\begin{claim}\label{cl:T0}
$|T_0|=k$.
\end{claim}

\begin{proof}[Proof of Claim~\ref{cl:T0}]
We first show that $|T_0|\ge k$.
Note that the initial factor degree of $x$ is $k$, and its
residual prescribed degree becomes $0$ after deleting the
vertex set $T_0$.
Thus $N_r(x)\subseteq T_0$, and
$|N_r(x)|=d_F(x)=k$.
Therefore $|T_0|\ge k$.

Next we prove that $|T_0|\le k$ by contradiction.
Suppose that $|T_0|>k$.
We construct an enlarged graph $G'$ by adding all edges from
$x$ to vertices in $T_0\setminus N_r(x)$.
These edges are not in $G$. Indeed, otherwise, if a vertex
$u\in T_0\setminus N_r(x)$ were adjacent to $x$, then $xu$
would be red, contrary to $u\notin N_r(x)$.

We now verify that $G'$ has a unique $k$-factor.
Since no edge is removed, $F$ remains a
$k$-factor of $G'$.
The $k$ vertices in $N_r(x)\subseteq T_0$ are unaffected by the
newly added edges, so they still have degree $k$
in $G'$.
Consequently, all their incident edges must
be included in any $k$-factor of $G'$.
In particular, the $k$ red edges between $N_r(x)$ and $x$ are
forced into every factor, giving $x$ its prescribed degree $k$.
Consequently, none of the newly added edges incident to $x$
can belong to any $k$-factor.
Thus the set of $k$-factors of $G'$ is identical to that of $G$,
which implies that the factor of $G'$ is unique.

On the other hand, one can see that $|E(G')|=|E(G)|+|T_0|-k>|E(G)|,$ contradicting the assumption that $G$ is a maximum-size graph with a unique $k$-factor.
Therefore $|T_0|=k$, completing the proof of the claim.
\end{proof}

Since $N_r(x)\subseteq T_0$ and both sets have
size $k$, we have $T_0=N_r(x)$.

We assert that $x$ is a universal vertex of $G$.
Repeated applications of
Lemma~\ref{lem:delete}(ii) show that $G-T_0$ is extremal for
the residual degree function
\[
 f_0(v)=d_{F-T_0}(v)=k-|N_r(v)\cap T_0|,
 \qquad \text{for }v\in V(G)\setminus T_0.
\]
Furthermore, $x$ is a free vertex in $G-T_0$.
By Lemma~\ref{lem:delete}(i), every free vertex in an extremal
graph is universal, so $x$ is adjacent to all other vertices
in $G-T_0$.
Since $x$ is also adjacent to every vertex in $T_0$
(as $T_0=N_r(x)$), it follows that $x$ is adjacent to every
other vertex in $G$.
Hence the assertion holds.

By the assertion of the universal vertex $x$ of $G$, $d_G(x)=n-1>k$, since $k\le n-2$.
For any other vertex $y\notin T_0\cup\{x\}$,
$xy$ is blue because $x$ is universal and
$N_r(x)=T_0$.
Together with the $k$ red edges incident to $y$, this gives
$d_G(y)\ge k+1$.
Thus, by Claim~\ref{cl:T0},
$D=T_0=N_r(x)$ and $|D|=k$.

The proof of Theorem~\ref{thm:structure} is complete.
\end{proof}

\end{document}